\documentclass{amsart}
\usepackage{epsf,amsfonts,amssymb}
\usepackage{epsfig}
\usepackage{amsthm}
\usepackage[mathscr]{eucal}
\def\t{{\large $\tau$}$\!_{1}$}    
\newcounter{fig}
\numberwithin{equation}{section}
\DeclareMathOperator{\Aut}{Aut}
\DeclareMathOperator{\diam}{diam}
\DeclareMathOperator{\dm}{d}

\newtheorem{thm}{Theorem}
\newtheorem{lem}{Lemma}
\newtheorem{cor}{Corollary}
\theoremstyle{remark}
\newtheorem{rem}{Remark}
\newtheorem{defn}{Definition}
\newtheorem{exmp}{Example}
\begin{document}

\title[Quasi-isometries between non-locally-finite graphs and structure trees]{Quasi-isometries between non-locally-finite graphs and structure trees}
\author[B.\ Kr\"on]{B.\ Kr\"on$^{\textstyle{\,\star}}$}
\begin{abstract}
We prove several criteria for quasi-isometry between non-locally-finite graphs and their structure trees. Results of M\"oller in \cite{moeller92ends2} for locally finite and transitive graphs are generalized. We also give a criterion which describes quasi-isometry by how edge-ends are split up by the cuts of a structure tree. 
\end{abstract}
\bibliographystyle{plain}
\thanks{$^{\textstyle{\star}}$ The author is supported by the START-project Y96-MAT of the Austrian Science Fund. Current address: Institut f\"ur Mathematik C, Technische Universit\"at Graz, Steyrergasse 30, 8010 - Graz, tel.: +43/316/873-4509, e-mail: kroen@finanz.math.tu-graz.ac.at\\
Mathematics Subject Classification 05C75, (05C25, 20B27)}

\maketitle

\section{Introduction}
Quasi-isometry on graphs is a weakened form of isomorphism. Graphs which are quasi-isometric to each other have the same global structure but may have local deviations which are uniformly bounded. The main property of quasi-isometries is that a set has finite diameter if and only if its image has finite diameter.\par
In \cite{gromov84infinite} and \cite{gromov93asymptotic}, Gromov used the concept of quasi-isometry in the context of structural properties of infinite groups.\par
By cutting a graph into pieces so that the resulting set of cuts is invariant under the action of the automorphism group we obtain a \emph{structure tree set}. The lines along which the graph is sliced into pieces can be regarded as level lines of a map. By identifying the ranges in that map that lie between these lines with new vertices and connecting pairs of them when they are separated by just one level line we obtain a tree called the \emph{structure tree}. These ranges can be empty but it will turn out that no two empty ranges can be adjacent.\par
If the structure tree is quasi-isometric to the graph then it describes the ramification structure of this graph. We could say the graph `looks like' the structure tree. In this article we give a detailed discussion of when a graph is quasi-isometric to a structure tree. The general criterion in Theorem \ref{general_qi-criterion} says that this is the case if and only if the diameter of the ranges mentioned above, together with their surrounding level lines, is finite.\par
There are two ways for a non-locally-finite graph to have infinite diameter: by having rays with infinite diameter or by having so-called \emph{star-balls}. The essence of Theorem \ref{general_Phi-criterion} is that a graph is quasi-isometric to a structure tree if and  only if a)~it has no star-balls and b)~any rays of infinite diameter are cut into pieces by the structure tree set.

\section{Structure trees}
Throughout this article let $X=(VX,EX)$ be a connected, undirected graph without loops or multiple edges. The set $VX$ of all vertices consists of $VX^{L}$, the set of vertices with finite degree, and the set of vertices with infinite degree which we denote by $VX^{\infty}$. A set $e$ of vertices in $VX$ is called connected, if any two vertices in $e$ can be connected by a path in $X$ that does not leave $e$. We write $e^{*}$ for the complement $VX\backslash e$ of $e$ and $\diam_{X}$ for the diameter with respect to the natural graph metric $\dm_{X}$ in $X$. The closed ball with center $x$ and radius $r$ is denoted by $B(x,r)$.\par
The \emph{vertex-boundary} $\theta e$ of $e$ is the set of vertices in $e^{*}$ which
are adjacent to a vertex in $e$. $I\theta e:= \theta e^{*}$ is called \emph{inner vertex-boundary} of $e$. The \emph{edge-boundary} $\delta e$ of $e$ is defined as the
set of edges connecting vertices in $e$ with vertices in $e^{*}$. A non-empty set of vertices $e$ is a \emph{cut} (or \emph{edge-cut}) if $\delta e$ is finite. For $n=|\delta e|$ we also call $e$ an \emph{$n$-cut}. If both $e$ and $e^{*}$ are connected, a cut $e$ is said to be \emph{tight}.

\begin{defn}
A set $E$ of cuts in $X$ is called a \emph{tree set}, if it satisfies the following three axioms.
\vskip6pt
\begin{enumerate}
\item[(S1)] For all pairs of cuts $e$ and $f$ in $E$, one of the following inclusions holds:
$$e\subset f,\quad e\subset f^{*},\quad e^{*}\subset f \mbox{\quad or \quad} e^{*}\subset f^{*}.$$\par
\item[(S2)] For any two cuts $e$ and $f$ in $E$ there exist only finitely many cuts $d$ in $E$ such that $e\subset d \subset f$.
\item[(S3)] Neither $\emptyset$ nor $VX$ is an element of $E$.
\end{enumerate}
The tree set $E$ is called \emph{undirected}, if also
\begin{enumerate}
\item[(S4)] $e$ is an element of $E$ if and only if $e^{*}$ is an element of $E$.
\end{enumerate}
\vskip6pt
An undirected tree set that consists only of tight $n$-cuts is called \emph{tight} tree set.\par
We call an edge-cut $e$ \emph{non-trivial}, if both $e$ and $e^{*}$ are infinite. Non-trivial and tight edge-cuts $e$ for which $\Aut(X) e \cup \Aut(X) e^{*}$ is a tree set are called \emph{structure cuts}. Such a tree set is called a \emph{structure tree set}.
\end{defn}

\begin{thm}
If a graph has a non-trivial cut then it also has a structure cut.
\end{thm}

This important theorem was originally stated by Dunwoody in \cite[Theorem 1.1]{dunwoody82cutting}. An improved version of the proof can be found in \cite{dicks89groups}.\par

\begin{defn}
Let $e$ and $f$ be cuts in a tree set $E$. We say $e$ \emph{points to} $f$ (notation $e\gg f$), if $f$ is a subset of $e$ and there is no third cut $d\in E$, such that $f\subset d \subset e$. The cuts \emph{point away} from each other (notation $e \rightleftharpoons f$) if $e^{*}\gg f$ and $f^{*}\gg e$ or $f=e$.
\end{defn}

The ramifications that can be described by a tree set can always be represented by a tree which is called a \emph{cut tree}. We want to give an axiomatic definition.

\begin{defn}
A \emph{cut tree} of a tree set $E$ is a connected directed tree $T=T(E)$, for which there exists a bijection $b:E\to ET$ with the following properties:
\vskip6pt
\begin{enumerate}
\item[(T1)] $b(e)=(u,v)$ is equivalent to $b(e^{*})=(v,u)$ and
\item[(T2)] $e\gg f$ is equivalent to $t(b(e))=o(b(f))$
\end{enumerate}
\vskip6pt
where $o(p)$ is the origin and $t(p)$ is the terminus of any directed edge $p$. If $E$ is a structure tree set then $T=T(E)$ is called a \emph{structure tree}.\par

\end{defn}
Note that (T1) an (T2) imply that $e\gg f$ is also equivalent to $t(b(f^{*}))=o(b(e^{*}))$. To avoid complicated notation we will not distinguish between a cut $e$ and the corresponding edge $b(e)$.

\begin{thm}
To every undirected tree set\/ $E$ there exists a cut tree\/ $T=T(E)$ which is unique up to isomorphism.
\end{thm}

The existence of cut trees for a given tree set was proved in \cite[Theorem 2.1]{dunwoody79accessibility}. Various examples for structure trees can be found in \cite[Section 2.3]{moeller95groups}.\par
The following lemma is a generalization of a statement of Dunwoody in \cite[2.3]{dunwoody82cutting}. Thomassen \cite[Proposition 4.1]{thomassen93vertex} found a surprisingly simple proof by induction.

\begin{lem}\label{thomassen}
For every given natural\/ $n$ and every edge\/ $p$ in a connected graph\/ $X$ there exist only finitely many tight n-cuts\/ $e$ such that\/ $p$ is element of the edge-boundary\/ $\delta e$.
\end{lem}

\begin{cor}\label{cor_of_lem_of_Th}
Every strictly decreasing sequence of tight n-cuts whose intersection is non-empty must be finite.
\end{cor}

\section{Edge-ends}

A \emph{ray} is a sequence $(x_{n})_{n\in \mathbb{N}}$ of pairwise distinct
vertices such that $x_{n}$ is adjacent to $x_{n+1}$ for all $n$. We write $RX$ for the set of all rays in $X$. A ray \emph{lies} in a set $e$ of vertices or is \emph{contained} in $e$, if $e$ contains all but finitely many elements of the ray. Sometimes we will use the terms \emph{contain} and \emph{lie} at the same time in the sense above as well as in the sense of set theoretic inclusion. A set $e$ of vertices \emph{separates} two sets of vertices or rays, if one of them lies in $e$ and the other lies in $e^{*}$.\par
Two rays are called \emph{edge-equivalent in the first sense} if they cannot be separated by \emph{edge-cuts}. It is easy to see that this relation is an equivalence relation. Its equivalence classes are called \emph{edge-ends of the first type}.\par
An end \emph{lies} in a set of vertices $e$ or is \emph{contained} in $e$, if all of its rays lie in $e$. The set of edge-ends of the first type that lie in $e$ is denoted by $\Omega_{1} e$. In fact, an edge-end of the first type $\omega$ lies in an edge-cut $e$ if and only if one of its rays lies in $e$. So $\omega$ either lies in $e$ or in $e^{*}$.
\begin{lem}
For a graph\/ $X$, the set
$$BX:=\{e\cup \Omega_{1} e \mid e \subset VX \mbox{\ and\ \,} |\delta e| <
\infty \}$$
is closed under finite intersection.
\end{lem}

For a proof see e.g.\ \cite[Lemma 8]{kroen00end}. $BX$ is a base of a topological space $(VX\cup\Omega_{1} X,$\t$X)$ whose topology \t$X$ is called \emph{edge-topology of the first type}. By Theorem 2 and Example 3 in \cite{kroen00end} we know that \t$X$ is compact but in general not even $T_{0}$. For graphs with countably finite degree this compactness can easily be deduced from results of Cartwright, Soardi and Woess in \cite{cartwright93martin}.\par
To obtain better properties of separation we will now extend the edge-equivalence to $RX\cup VX^{\infty}$. This strategy was first adopted in the article \cite{cartwright93martin} mentioned above.\par
Two elements of $RX\cup VX^{\infty}$ are called \emph{edge-equivalent in the second sense} if for every edge-cut $e$ either both lie in $e$ or both lie in $e^{* }$. Again it is easy to see that this relation is an equivalence-relation. We call its equivalence-classes \emph{edge-ends} or \emph{edge-ends of the second type}. The terms \emph{to lie in} and \emph{separate} are used in the same sense as above. The set of all ends lying in some set of vertices $e$ is denoted by $\Omega e$. We usually write $\Omega X$ instead of $\Omega VX$. A finite set of vertices $e$ \emph{separates} two ends if they lie in different connected components of $e^{*}$.\par
In every edge-end of the second type containing the same ray, there lies an edge-end of the first type. But note that there also may exist ends of the second type consisting only of vertices. In \cite{cartwright93martin} these ends are called \emph{improper ends}.\par
By the same construction as in the first case we now obtain the \emph{edge-topology of the second type}. It is normal, Lindel\"off and totally disconnected, see \cite{kroen00end}. Compactness can be deduced from the compactness of the edge-topology of the first type.

\section{Vertex and end structure mapping}

To describe the connections between a graph $X$ and its structure trees $T=T(E)$ we now want to define functions $\phi: VX\to VT$ and $\Phi: \Omega X\to \Omega T$. Another construction of $\phi$ can be found in \cite{dicks89groups}. M\"oller gave a construction for the function $\Phi$ in \cite[Proposition 1]{moeller92ends2} for pairs of quasi-isometric graphs in the locally finite case. In \cite[Section 7]{kroen00end} the author studied the connections between quasi-isometries and a similar function on another end compactification, the so-called \emph{metric end compactification} of non-locally-finite graphs.

\begin{defn}
Let $T=T(E)$ be a cut tree of a graph $X$. A cut $e$ in $E$ \emph{points at} some vertex $x$ in $VX$ (notation: $e\to x$), if $x$ is an element of $e$ and there is no other cut which contains $x$ and is a subset of $e$.
\end{defn}

\begin{lem}\label{unique_terminal_vertex}
For every\/ $x\in VX$ there exists a cut\/ $e$ in a tight tree set\/ $E$ such that\/ $e\to x$. The cuts that point at\/ $x$, seen as edges in\/ $T$, have all the same terminal point.
\end{lem}

\begin{proof}
Since a tight tree set is undirected there must exist a cut $e_{1}$ in $E$ that contains $x$. By Corollary 1 every sequence $(e_{n})_{n\ge 1}$ of cuts in $E$ containing $x$ must be finite. The last cut in such a sequence of maximal length must point at $x$.\par
Suppose that there are two cuts $e_{1}$ and $e_{2}$ in $E$, such that $e_{1}\to x$, $e_{2}\to x$ and $t(e_{1})\ne t(e_{2})$. By Axiom S1 of the definition of a tree set we distinguish between four cases. $e_{1}\subset e_{2}$, $e_{2}\subset e_{1}$ and $e_{1}\cap e_{2}=\emptyset$ would immediately imply a contradiction. If $e_{1}\cup e_{2}=VX$ then there must exist some $f\in E$ such that $e_{1}^{*}\subset f \subset e_{2}$ since $t(e_{1})\ne t(e_{2})$. If $x\in f$ then $e_{2}$ does not point at $x$; if $x$ does not lie in $f$ then $e_{1}$ does not point at $x$.
\end{proof}

\begin{defn}
Let $\phi$ be the function $VX\to VT$ such that $\phi(x):=t(e)$ for some cut $e\in E$ which points at $x$. We call $\phi$ the \emph{vertex structure mapping} with respect to $T$.\par
In a similar way we now want to construct the \emph{end structure mapping} $\Phi: \Omega X\to VT\cup \Omega T$. For each $\omega\in \Omega X$ we have two cases.
\begin{enumerate}
\item For every cut $e\in E$ containing $\omega$ there exists another cut $f\in E$ in which $\omega$ is contained and for which $e \gg f$. In this case any decreasing sequence of cuts in $E$ containing $\omega$ defines a unique end $\varepsilon\in \Omega T$ and we set $\Phi(\omega):=\varepsilon$.
\item The end $\omega$ lies in some $e\in E$ but in no further cut in $E$ which is contained in $e$. We set $\Phi(\omega):=t(e)$. The uniqueness of $t(e)$ can be seen by the same arguments that we used in the second part of the proof of Lemma 3.
\end{enumerate}
\end{defn}

By Lemma \ref{unique_terminal_vertex} the definition of the vertex structure mapping is independent of the choice of $e$.

\begin{lem}[{\cite[Lemma 2]{moeller92ends1}}]\label{qi_the_Phi_bij_on...}
For a tight cut tree\/ $T(E)$ the restriction of\/ $\Phi$ on\/ $\Phi^{-1}(\Omega T)$ is bijective.
\end{lem}

\section{The action of Aut(X) on a structure tree}\label{act_of_Aut_on}

For a vertex $x$ in $VX$ and a tree set $E$ we define $N(x):=\{e \in E\mid e\to x\}$. If $g$ is an automorphism of $X$ then we have
$$g N(x)=g\{e \in E\mid e\to x\}=\{ge \in E\mid e\to x\}.$$
Since the cut $e$ points to $x$ if and only if $ge$ points to $gx$, this set is equal to
$$\{ge \in E\mid ge\to gx\}=\{f \in E\mid f\to gx\}=N(gx).$$
The images $g\phi^{-1}\phi(x)$ and $\phi^{-1}\phi g(x)$ are the sets of all vertices pointed at by cuts in $gN(x)$ or $N(gx)$, respectively. We now define a function
$$\bar{g}^{T}:\phi(VX)\to\phi(VX),\ v\mapsto \phi g\phi^{-1}(v).$$
By the above considerations we obtain
$$\bar{g}^{T}\phi(x)=\phi g \phi^{-1}\phi(x)=\phi\phi^{-1}\phi g(x)=\phi g(x).$$
For all $x\in VX$ and $v\in \phi(VX)$ we now have the following formulas

\begin{align}\label{3_equs}
g \phi^{-1}\phi(x)&=\phi^{-1}\phi g(x),\\
\phi^{-1}\bar{g}^{T}(v)&=g\phi^{-1}(v),\notag \\
\bar{g}^{T}\phi(x)&=\phi g(x).\notag
\end{align}

Thus $\bar{g}^{T}$ is a well defined function which is induced in a natural way by the automorphism $g$ of $X$.\par
If we assume that $\phi(VX)$ does not cover the whole set of vertices $VT$ then  $\phi(VX)$ is one of the bipartite blocks in $T$. For any two $\phi$-images $\phi(x)$ and $\phi(y)$ at distance 2 in $T$ we can find cuts $e$ and $f$ in $ET$ such that $e\to x$, $f\to y$ and $e\leftrightharpoons f$. There also exist cuts with this property for the vertices $g(x)$ and $g(y)$ in $VX$. By \ref{3_equs} we now obtain
$$\dm_{T}(\phi(x),\phi(y))=\dm_{T}(\phi g(x),\phi g(y))=\dm_{T}(\bar{g}^{T}\phi(x),\bar{g}^{T}\phi(y))=2.$$
Since $T$ is a tree this implies the following

\begin{lem}\label{d_T_equation_for_phi-images}
For all pairs of vertices\/ $x$ and\/ $y$ in\/ $VX$
$$\dm_{T}(\phi(x),\phi(y))=\dm_{T}(\phi g(x),\phi g(y))=\dm_{T}(\bar{g}^{T}\phi(x),\bar{g}^{T}\phi(y)).$$
\end{lem}

The function $\bar{g}^{T}$ now can easily be extended to a bijective isometry on the whole set of vertices $VT$. This automorphism of $T$ is denoted by $g^{T}$.\par
If $\phi(VX)=VT$ we can see by the same arguments that $\bar{g}^{T}$ itself is already an automorphism of $T$. In this case we define $g^{T}=\bar{g}^{T}$.\par
The set
$$\Aut^{T}(X):=\{g^{T}\mid g\in \Aut(X)\}$$
of these functions acts transitively at least on the bipartite blocks of $T$. It acts transitively on the whole structure tree if and only if there exist cuts in $\Aut(X)e$ as well as in $\Aut(X)e^{*}$ that both point at some vertex $x$ in $VX$.\par
To see that in the general case the function
$$L:\Aut(X)\to \Aut(T): g\mapsto {g}^{T}$$
is neither surjective nor injective we give the following example.

\begin{exmp}
At each vertex of the cycle $C_{4}=(v_{1},v_{2},v_{3},v_{4})$ of length 4 we fix a pair of hanging edges. The unions of the vertices in these pairs of hanging edges and their complements constitute a tree set $E$ with eight elements. $T(E)$ is isomorphic to the star $K_{1,4}$. Its vertex of degree 4 is denoted by $v$. Instead of pairs of oppositely oriented edges in $ET$ we draw undirected edges. See Figure \ref{fig_stru_tree}.
\begin{center}
\setlength{\unitlength}{1mm}
\begin{picture}(140,40)
\refstepcounter{fig}
\label{fig_stru_tree}
\put(45,0){\circle*{1}}
\put(55,0){\circle*{1}}
\put(35,10){\circle*{1}}
\put(45,10){\circle*{1}}
\put(55,10){\circle*{1}}
\put(65,10){\circle*{1}}
\put(35,20){\circle*{1}}
\put(45,20){\circle*{1}}
\put(55,20){\circle*{1}}
\put(65,20){\circle*{1}}
\put(55,30){\circle*{1}}
\put(45,30){\circle*{1}}
\put(100,5){\circle*{1}}
\put(100,25){\circle*{1}}
\put(120,5){\circle*{1}}
\put(120,25){\circle*{1}}
\put(110,15){\circle*{1}}
\put(35,10){\line(1,0){30}}
\put(35,20){\line(1,0){30}}
\put(45,0){\line(0,1){30}}
\put(55,0){\line(0,1){30}}
\put(100,5){\line(1,1){20}}
\put(100,25){\line(1,-1){20}}
\put(30,30){$X$}
\put(86,30){$T(E)$}
\put(109,18){$v$}
\put(40,6.5){$v_{3}$}
\put(40,22){$v_{2}$}
\put(57,6.5){$v_{4}$}
\put(57,22){$v_{1}$}
\linethickness{0.25mm}
\bezier{20}(32,18)(53,12)(47,33)
\bezier{20}(32,12)(53,18)(47,-3)
\bezier{20}(68,18)(47,12)(53,33)
\bezier{20}(68,12)(47,18)(53,-3)
\put(70,-10){Figure \ref{fig_stru_tree}}
\end{picture}
\end{center}
\vskip1.3cm

Each permutation of $VT\backslash \{v\}$ corresponds to an automorphism of the tree $T$, whereas automorphisms of $X$ must respect the structure of the cycle $C_{4}$. Thus $L$ is not surjective.\par
Automorphisms that have the same action on $C_{4}$, but different action on the vertices of degree one in $VX$ are all mapped to the same automorphism of $T$. This means that the operator $L$ is not injective.
\end{exmp}

\begin{lem}
For vertices\/ $v$ and\/ $w$ in\/ $\phi(VX)$ and an automorphism\/ $g$ in\/ $\Aut (X)$ we have
$$g^{T}(v)=w \iff g\phi^{-1}(v)=\phi^{-1}(w).$$
\end{lem}

\begin{lem}
For vertices $v$ and $w$ in $\phi(VX)$ and an automorphism $g$ in $\Aut (X)$ we have
$$g^{T}(v)=w \iff g\phi^{-1}(v)=\phi^{-1}(w).$$
\end{lem}

\begin{proof}
By (5.1) $g\phi^{-1}(v)=\phi^{-1}(w)$ is equivalent to $\phi^{-1}g^{T}(v)=\phi^{-1}(w)$.
\end{proof}

\begin{lem}\label{same_diam_of_phi^{-1}_in_blocks}
For any two vertices\/ $v$ and\/ $w$ lying in some bipartite block of\/ $T$ we have
$$\diam_{X}\phi^{-1}(v)=\diam_{X}\phi^{-1}(w).$$
\end{lem}

\begin{proof}
If these $\phi$ pre-images are non-empty there exists a $g^{T}\!\in\Aut^{T}(X)$ for which $g^{T}(v)\!=\nolinebreak\!w$. The statement now is a consequence of $g\phi^{-1}(v)=\phi^{-1}(w)$.
\end{proof}

\begin{defn}
For a vertex $v$ in a structure tree $T$ we write $N(v)$ for the set of all cuts $e$ in $ET$ with $t(e)=v$. We write $N(v)^{*}$ for the set of all cuts $f$ for which $f^{*}$ is an element of $N(v)$. The set $$R(v):=\phi^{-1}(v)\cup \{x\mid x\in I\theta e, e\in N^{*}(v)\}=\phi^{-1}(v)\cup \{x\mid x\in \theta e, e\in N(v)\}$$
is called the \emph{region of $v$}.
\end{defn}

If $\phi^{-1}(v)$ is nonempty and $x\in \phi^{-1}(v)$ then $N(x)=N(v)$.

\begin{lem}\label{same_diam_of_R_in_blocks}
For any two vertices\/ $v$ and\/ $w$ in the same bipartite block of\/ $T$,
$$\diam_{X}R(v)=\diam_{X}R(w).$$
\end{lem}
\begin{proof}
We can find an authomorphism $g\in \Aut(X)$ such that $g(v)=w$ and $R(w)=R(g^{T}(v))=gR(v)$.
\end{proof}

\section{A general criterion for quasi-isometry of $\phi$}

\begin{defn} \label{qi-defn}
Two graphs $X$ and $Y$ are called \emph{quasi-isometric with respect to the functions} $\phi: VX\to VY$ \emph{and} $\psi: VY\to VX$ if there exist constants $a$, $b$, $c$ and $d$ such that for all vertices $x$,
$x_{1}$ and $x_{2}$ in $VX$ and vertices $y$, $y_{1}$ and $y_{2}$ in $VY$, the following conditions hold
\begin{enumerate}
\item[(Q1)] $\dm_{Y}(\phi(x_{1}),\phi(x_{2}))\le a\cdot\dm_{X}(x_{1}, x_{2})$ \quad (boundedness of $\phi$)
\item[(Q2)] $\dm_{X}(\psi(y_{1}),\psi(y_{2}))\le b\cdot\dm_{Y}(y_{1}, y_{2})$ \quad (boundedness of $\psi$)
\item[(Q3)] $\dm_{X}(\psi\phi(x),x)\le c$ \quad (quasi-injectivity of $\phi$)
\item[(Q4)] $\dm_{Y}(\phi\psi(y),y)\le d$ \quad (quasi-surjectivity of $\psi$)
\end{enumerate}
We call $\phi$ and $\psi$ \emph{quasi-isometries}. They are said to be \emph{quasi-inverse} to each other.
\end{defn}
For general metric spaces the definition of quasi-isometry includes further additive constants in the Axioms (Q1) and (Q2). In case the positive values of the metric are greater than some positive real number these additive constants are not needed.\par
Quasi-isometries may change structures as long as the differences can be bounded uniformly. In other words we could say that they preserve the global structure of graphs when we consider graphs as discrete metric spaces only.\par
Quasi-isometry is an equivalence relation on the family of all graphs. Various examples of quasi-isometric graphs can be found in \cite[Example 6]{kroen00end}.\par
The following lemma describes the most important of the basic properties of quasi-isometries.

\begin{lem}\label{qi_and_finite_diam}
Let\/ $\phi: VX\to VY$ be a quasi-isometry and\/ $A$ a subset of\/ $VX$. Then
$$\diam_{X}A<\infty \Leftrightarrow \diam_{Y}\phi(A)<\infty.$$
\end{lem}

The following extends a result of M\"oller \cite[Lemma 1]{moeller92ends2} from locally finite graphs to arbitrary graphs.

\begin{thm}\label{general_qi-criterion}
A connected graph is quasi-isometric to a structure tree by the vertex structure mapping if and only if the regions of origin and terminus of an edge (equivalently: all edges) in the tree are bounded.
\end{thm}

\begin{proof}
Let $e$ be a cut in a structure tree set $E$ of the structure tree $T$. By Lemma \ref{thomassen} there is a natural number $k$ such that for every edge $(x,y)\in \delta e$ there exist no more than $k$ cuts in $E$ which contain $x$ but do not contain $y$. This constant is the same for all cuts in $E$. By this argument and by $\diam_{T}\phi \phi^{-1}(v)=0$ we obtain $\diam_{T} \phi R(v)\le 2k$ for every vertex $v\in VT$. Thus, by Lemma \ref{qi_and_finite_diam}, $\phi$ cannot be a quasi-isometry if the region of any vertex in $VT$ has infinite diameter.\par
Now we assume that there exists a cut $e$ in $E$ such that $R(o(e))$ and $R(t(e))$ both have finite diameter in $X$. By Lemma \ref{same_diam_of_R_in_blocks} this implies that all regions of vertices in $VT$ have finite diameter in $X$. To prove the theorem we now have to show that $\phi$ is a quasi-isometry. We now construct a function $\psi:VT\to VX$ which will turn out to be a quasi-inverse of $\phi$.\par
If a vertex $v$ in $VT$ is not contained in $\phi(VX)$ then let $r(v)$ be an arbitrary vertex in $VT$ which is adjacent to $v$. Otherwise we set $r(v):=v$. If $\psi(v)$ is an arbitrary element of $\phi^{-1}(r(v))$ then we have $\phi\psi (v)=r(v)$. We will now prove that $\phi$ and $\psi$ are quasi-isometries that are quasi-inverse to each other by checking the four axioms in Definition \ref{qi-defn}.\par
(Q1) \quad If there is no cut $e$ in $E$ which separates two adjacent vertices $x$ and $y$ then $\phi(x)=\phi(y)$. Otherwise there must be an automorphism $g\in \Aut (X)$ such that $g(\{x,y\})$ is in the edge-boundary $\delta e$. By Lemma \ref{d_T_equation_for_phi-images} and since $\delta e$ is finite the distance $\dm_{T}(\phi (x),\phi (y))$ can have only finitely many values for adjacent vertices $x$ and $y$. Thus the set
$$\{\dm_{T}(\phi (x),\phi (y))\mid \dm_{X}(x,y)=1\}$$
has a maximal element $a$. For any two vertices $x$ and $y$ in $VX$ there is a path \linebreak $\{x= z_{0},z_{1},\dots,z_{n}= y\}$ of length $\dm_{X}(x,y)$ and thus we obtain
$$\dm_{T}(\phi (x),\phi (y))\le
\dm_{T}(\phi (x=z_{0}),\phi (z_{1}))+\dots +\dm_{T}(\phi (z_{n-1}),\phi (z_{n}=y))\le
a \cdot \dm_{X}(x,y).$$

(Q2) \underline{case 1.} \quad $\phi(VX)\ne VT$\par
Let $v_{1}$ and $v_{2}$ be two vertices in $\phi(VX)$ with $\dm_{T}(v_{1},v_{2})=2$ and let $w$ be the vertex in $VT\backslash \phi(VX)$ which is adjacent to $v_{1}$ and $v_{2}$. Then
$$\diam_{X}\phi^{-1}(v_{1})+\diam_{X}\phi^{-1}(v_{2})+\dm_{X}(\phi^{-1} (v_{1}),\phi^{-1} (v_{2}))\le$$
$$\diam_{X}\phi^{-1}(v_{1})+\diam_{X}\phi^{-1}(v_{2})+\diam_{X}R(w).$$
By Lemma \ref{same_diam_of_phi^{-1}_in_blocks} and Lemma \ref{same_diam_of_R_in_blocks} the latter sum does not depend on the choice of the vertices. Thus it is a constant which we denote by $2b$. For two vertices $x$ and $y$ in $VX$ and a path $(v_{0}=\phi(x),v_{1},v_{2},\dots,v_{2k}=\phi(y))$ of length $\dm_{T}(\phi(x),\phi(y))=2k$ such that $v_{2i}\in \phi(VX)$ for $0\le i\le k$ we finally have
\begin{align}
\dm_{X}(x,y)&\le \sum_{i=0}^{k-1} \left( \diam_{X}\phi^{-1}(v_{2i})+\diam_{X}\phi^{-1}(v_{2(i+1)})+
\dm_{X}( \phi^{-1}(v_{2i}),\phi^{-1}(v_{2(i+1)})) \right)\nonumber\\
&\le k\cdot 2b=b\cdot \dm_{T}(\phi(x),\phi(y)).\nonumber
\end{align}
\par
\hskip9.5mm \underline{case 2.} \quad $\phi(VX)= VT$\par
In this case the proof of Case 1 works analogously by using the inequality
$$\diam_{X}\phi^{-1}(v)+\diam_{X}\phi^{-1}(w)+\dm_{X}(\phi^{-1} (v),\phi^{-1} (w))$$
$$\le\diam_{X}R(v)+\diam_{X}R(w)$$
for any adjacent vertices $v$ and $w$ in $VT$.

(Q3) \quad A vertex $x$ in $VX$ and the vertex $\psi\phi(x)$ always lie in the same $\phi$-pre-image of some vertex in $VT$. By defining
$$c:=\max\{\diam_{X}\phi^{-1}(t(e)),\diam_{X}\phi^{-1}(t(e^{*}))\}$$
for any $e\in E$ we obtain $\dm_{X} (\psi\phi(x),x)\le c$ for all vertices $x$ in $VX$.

(Q4) \quad By $\phi\psi(v)=r(v)$ we obtain $\dm_{T}(\phi\psi(v),v)\le 1$ for all vertices $v$ of the structure tree $T$.
\end{proof}

Since we did not use the assumption $\diam_{X}R(v)<\infty$ in verifying Axiom (Q1) we have proved the following lemma.

\begin{lem}\label{general_Q1}
For a structure tree $T$ there exists a constant $a$ such that
$$\dm_{T}(\phi (x),\phi (y))\le a \cdot \dm_{X}(x,y)$$
for all vertices $x$ and $y$ in $VX$.
\end{lem}

\begin{lem}\label{diam_phi^-1_and_B}
If $\phi^{-1}(v)$ is non-empty for some vertex $v$ of a structure $T=T(E)$ of a graph $X$, then $R(v)$ has finite diameter if and only if $\phi^{-1}(v)$ has finite diameter.
\end{lem}

\begin{proof}
Since $\phi^{-1}(v)$ is a subset of $R(v)$ we need only prove that $\diam_{X}R(v)$ is finite if $\phi^{-1}(v)$ has finite diameter. By the definition of a structure tree set the stabilizer $\Aut_{\phi^{-1}(v)}(X)$ of the set $\phi^{-1}(v)$ has at most two orbits $O_{1}$ and $O_{2}$ on $N(v)^{*}$. For every cut $e$ in  $O_{1}$ the set
$$\{\dm_{X}(x,\phi^{-1}(v))\mid x\in I\theta e\}$$
has a maximal element. The same holds for the orbit $O_{2}$. The larger of these two maxima is the maximal distance between a vertex in $R(v)$ and a vertex in $\phi^{-1}(v)$. Since $\phi^{-1}(v)$ has finite diameter the same must hold for $R(v)$.
\end{proof}

The following example shows that for a connected graph and its structure tree to be quasi-isometric, it is not enough to have finite diameters of the $\phi$-pre-images.

\begin{exmp}
For the two-sided infinite line $L$ with $VL=\{x_{k}\mid k\in \mathbb{Z} \}$, the set $E_{L}=\{ \{x_{k}\} \cup \{x_{k}\}^{*}\mid k\in \mathbb{Z} \}$ is a structure set. The corresponding structure tree $T$ looks like a star with one vertex $v$ of infinite degree and infinitely many vertices of degree one. The $\phi$-pre-images of the vertices with degree one consist of one vertex whereas $\phi^{-1}(v)$ is empty. All these pre-images have finite diameter but $R(v)$ equals $VL$, and therefore it has infinite diameter. Thus $L$ and $T$ are not quasi-isometric to each other.
\end{exmp}

\refstepcounter{fig}
\label{comb}

\section{Uniform ramification}

When we want to find criteria for quasi-isometry between graphs and their structure trees by the end structure mapping $\Phi$ we have to take into consideration the fact that infinite diameters in non-locally-finite graphs do not necessarily occur in connection with rays. For an example of a graph with infinite diameter which does not contain any ray see \mbox{\cite[Example 1]{kroen00end}.}

\begin{defn}
For a set of vertices $B$ in $VX$ we write $\mathcal{C}(B)$ for the set of all connected components of $B^{*}$ and $\mathcal{C}_{0}(B)$ for the set of all components in $B^{*}$ with finite diameter. A \emph{star ball} in a graph $X$ is a ball $S$ for which
$$\sup\{\diam_{X} C\mid C \in \mathcal{C}_{0}(S)\}=\infty.$$
A graph has \emph{uniform ramification} if it is connected, has infinite diameter and does not contain a star ball. A ray that does not contain infinite sets of vertices of finite diameter is called \emph{metric}.
\end{defn}

\begin{lem}\label{ball_contains_sb_then_sb}
Every ball\/ $B$ which contains some star ball\/ $S$ is a star ball.
\end{lem}

\begin{proof}
Let $z$ be the centre of $S$. In $\bigcup \mathcal{C}_{0}(S)$ there is a sequence $(x_{n})_{n\in\mathbb{N}}$ of vertices such that $\dm_{X}(x_{n},z)=n$. Since only finitely many vertices of the sequence lie in $B$ and for all other elements of the sequence there exist components in $\mathcal{C}_{0}(B)$ which contain them, $B$ must be a star ball, too.
\end{proof}

\begin{lem}\label{exist_of_com_with_inf_diam}
The complement\/ $B_{1}^{*}$ of every ball\/ $B_{1}$ in a graph with uniform ramification contains a connected component\/ $C$ with infinite diameter. For every ball\/ $B_{2}$ containing\/ $B_{1}$ there is also a connected component in\/ $C\backslash B_{2}$ which has an infinite diameter.
\end{lem}

\begin{proof}
A graph that ramifies uniformly has infinite diameter. If $B_{1}^{*}$ consisted only of components with finite diameter then $B_{1}$ would be a star ball.\par
If there were only connected components with finite diameter in $C\backslash B_{2}$ then $B_{2}$ would be a star ball since $C\backslash B_{2}$ has infinite diameter.
\end{proof}

\begin{lem}\label{unif_ram_then_ray_in_comp}
Let $X$ be a graph with uniform ramification. Then every component of infinite diameter in the complement of a ball contains a metric ray.
\end{lem}

\begin{proof}
For some $x_{0}$ in $VX$ let $C_{1}$ be a component of infinite diameter in $\mathcal{C}(\{x_{0}\})$. By induction we now choose a sequence $(C_{n})_{n\in \mathbb{N}}$ of components having infinite diameter such that
$$C_{n}\in \mathcal{C}(B(x_{0},n)) \mbox{\rm \ and\ } C_{n+1}\subset C_{n}.$$
The existence of such a sequence is a consequence of Lemma \ref{exist_of_com_with_inf_diam}. Let $x_{1}$ be an arbitrary vertex in the inner vertex-boundary $I\theta C_{1}$. Since $C_{1}$ is connected we can find a path from $x_{1}$ to some vertex $x_{2}$ in $I\theta C_{2}$ that does not leave $C_{1}\backslash C_{2}$. Again by induction we obtain a sequence of paths whose union is a metric ray.
\end{proof}

The following lemma, which can also be found in \cite{kroen00end} in a slightly modified version, characterizes graphs with infinite diameter. Its second part is a corollary of Lemma \ref{unif_ram_then_ray_in_comp}. The first part is well known and also easy to be proved.\par
In \cite{halin66graphen} Halin characterises rayless graphs in a similar way, but without taking their metric into consideration.

\begin{lem}\label{diam_of_graphs}
\label{starballs,_rays_and_diam}
\mbox{}\vspace{-4pt}
\begin{enumerate}
\item  A locally finite graph has infinite diameter if and only if it contains a ray.
\item The diameter of a non-locally-finite graph is infinite if and only if it contains a metric ray or a star ball.
\end{enumerate}
\end{lem}

\section{Almost transitive graphs}

Usually almost transitive graphs are defined as graphs with only finitely many orbits of vertices under the action of the automorphism group. Our definition is based on the natural metric of graphs. It is equivalent to the definition above in the locally finite case but includes a bigger class of graphs in the non-locally-finite case.

\begin{defn}
The automorphism group $\Aut(X)$ of a graph $X$ acts \emph{almost transitively} on $VX$ if there exists a vertex $x_{0}\in VX$ and a constant $r(x_{0})\in \mathbb{N}$ such that $\dm_{X}(\Aut(X) x_{0},x)\le r(x_{0})$ for all vertices $x$ in $VX$. The graph $X$ is called \emph{almost transitive} if $\Aut(X)$ acts almost transitively on $VX$. We call a ball $B$ \emph{covering ball} of $X$ if
$$\bigcup_{g\in \Aut(X)} gB=VX.$$
\end{defn}

\begin{rem}
A graph is almost transitive if and only if it contains a covering ball.
\end{rem}

\begin{lem}\label{if_qi_then_at}
A graph $X$ which is quasi-isometric to a structure tree by the vertex structure mapping $\phi$ is almost transitive.
\end{lem}

\begin{proof}
The set of automorphisms $\Aut^{T}(X)$ acts transitively on both bipartite blocks of $T$. Every ball $B$ in $T$ with radius at least 2 is a covering ball of $T$. Since $\phi$ is a quasi-isometry the pre-image $\phi^{-1}(B)$, by Lemma \ref{qi_and_finite_diam}, has a finite diameter. Every ball in $X$ containing $\phi^{-1}(B)$ is a covering ball.
\end{proof}

\begin{lem}\label{al_trans_then_ram_unif}
Every connected almost transitive graph $X$ is uniformly ramifying.
\end{lem}

\begin{proof}
We assume that there is a star ball $B(z,n)$. Let $r$ be a radius such that $B(z,r)$ is a covering ball. In $\mathcal{C}_{0}(B(z,n))$ there is a component $C$ containing a vertex $\bar{y}$ whose distance to $B(z,n)$ is greater then $2n+r$. Let $y$ denote a vertex of the $z$-orbit such that $\bar{y}$ is an element of the covering ball $B(y,r)$. $B(y,n)$ is again a star ball which is now contained in $C$. Since $C^{*}$ is connected it is completely contained in one of the components of $B(y,n)^{*}$. All other components of $B(y,n)^{*}$ are contained in $C$. Thus $B(y,n)$ cannot be a star ball.
\end{proof}

\section{The general $\Phi$-criterion}
The arguments in the proof of the following lemma are similar to those of Theorem 6 in \cite{kroen00end}.
\begin{lem}\label{if_qi_then_rays...}\par
Let\/ $X$ be a graph which is quasi-isometric to a structure tree\/ $T$ by the vertex structure mapping\/ $\phi$ and let\/ $\omega$ be an end in\/ $\Omega X$.\par
\begin{enumerate}
\item If\/ $\omega$ contains a ray with infinite diameter then\/ $\Phi(\omega)\in \Omega T$.
\item If\/ $\omega$ contains a vertex then\/ $\Phi(\omega)\in VT$.
\end{enumerate}
\end{lem}
\begin{proof}
By connecting the $\phi$-images of adjacent vertices of a ray $L_{1}$ of infinite diameter in $\omega$ by geodesic paths whose lengths are at most the constant $a$ in Axiom (Q1) of quasi-isometry, we obtain a path $P$ in $T$. By Lemma \ref{qi_and_finite_diam} its diameter is infinite. Again by Lemma \ref{qi_and_finite_diam} all $\dm_{T}$-balls in $P$, as a subgraph of $T$, contain at most finitely many $\phi$-images of vertices in $L_{1}$. Since we have constructed $P$ only with paths of length at most $a$, the subgraph $P$ of $T$ must be locally finite. By Lemma \ref{diam_of_graphs}, $P$ must contain some ray $L_{2}$. The end of $L_{2}$ has to be the $\Phi$-image of $\omega$.\par
By Corollary \ref{cor_of_lem_of_Th} there is no infinite sequence of cuts in a structure cut set with nonempty intersection. Thus an end containing a vertex in $X$ cannot be mapped by $\Phi$ onto an end of $T$.
\end{proof}

\begin{defn}
An end is called \emph{thick} if it contains infinitely many disjoint rays. An end which is not thick is called \emph{thin}. Denote by $\Theta X$ is the set of thick ends in a graph $X$ and by $\Delta X$ the set of thin ends. An end that does only contain vertices and rays of finite diameter is called a \emph{point end}. A \emph{mixed end} contains a ray with infinite diameter and a vertex of infinite degree. All other ends are called \emph{proper ends}. The set of point ends is denoted by $\Omega_{0}X$, the set of mixed ends by $\Omega_{1}X$ and the set of proper ends by $\Omega_{2}X$. Furthermore we define\par

$$\Delta_{2}X:=\Delta X \cap \Omega_{2}X,$$
$$\Theta_{2}X:=\Theta X \cap \Omega_{2}X.$$
\end{defn}
For the following observations it will not be necessary to distinguish between thick and thin ends in $\Omega_{0}X$ and $\Omega_{1}X$.
\begin{exmp}
A graph with one thin mixed end.

\begin{center}
\setlength{\unitlength}{1mm}
\begin{picture}(45,34)
\refstepcounter{fig}
\label{thin_mixed_end}
\multiput(0,10)(10,0){5}{\circle*{1}}
\multiput(10,15)(10,0){4}{\circle*{1}}
\multiput(20,20)(10,0){3}{\circle*{1}}
\multiput(30,25)(10,0){2}{\circle*{1}}
\put(40,30){\circle*{1}}
\put(25,0){\circle*{1}}
\put(25,0){\line(-5,2){25}}
\put(25,0){\line(-3,2){15}}
\put(25,0){\line(-1,2){5}}
\put(25,0){\line(1,2){5}}
\put(25,0){\line(3,2){15}}
\put(0,10){\line(1,0){40}}
\put(0,10){\line(2,1){40}}
\put(10,10){\line(0,1){5}}
\put(20,10){\line(0,1){10}}
\put(30,10){\line(0,1){15}}
\put(40,10){\line(0,1){20}}
\linethickness{0.25mm}
\bezier{6}(40,30)(43,31.5)(46,33)
\bezier{19}(25,0)(34.375,3.75)(43.75,7.5)
\bezier{6}(40,10)(43,10)(46,10)
\put(18,-8){Figure \ref{thin_mixed_end}}
\end{picture}
\end{center}
\end{exmp}
\vskip8mm
\begin{lem}
Proper ends only consist of rays with infinite diameter.
\end{lem}

\begin{proof}
We have to prove that a ray $L$ which is not equivalent to any vertex has infinite diameter. Let $C_{0}$ be a cut containing $L$. Every vertex $x$ in $I\theta C_{0}$ can be separated from $L$ by a cut $D_{x}$. The intersection
$$C_{1}:=C_{0}\backslash\bigcup \{D_{x}\mid x\in I\theta C_{0}\}$$
is again an edge-cut containing $L$. By induction we obtain a strictly decreasing sequence $(C_{n})_{n\in \mathbb{N}}$ with empty intersection such that $L$ lies in all the cuts $C_{n}$ and the inner vertex-boundaries of these cuts are pairwise disjoint. The distance of a vertex in $\theta C_{0}$ to any vertex in $C_{n}$ is greater than $n$. Since $L$ lies in all cuts $C_{n}$ it must have an infinite diameter.
\end{proof}

For the end structure mapping $\Phi$ we define two properties:
\begin{enumerate}
\item[(P1)] $$\Phi^{-1}(VT)=\Omega_{0}X$$
\item[(P2)] $$\Phi^{-1}(VT)=\Omega_{0}X,\quad \Omega_{1}X=\emptyset,\quad \Theta_{2}X=\emptyset, \mbox{\quad and\quad } \Phi(\Delta_{2}X)=\Omega T$$
\end{enumerate}
\begin{thm}\label{general_Phi-criterion}
For a graph\/ $X$ with a structure tree\/ $T(E)$ the following statements are \mbox{equivalent.}
\begin{enumerate}
\item $X$ is quasi-isometric to\/ $T$ by the vertex structure mapping $\phi$.
\item $X$ is uniformly ramifying and has property (P1).
\item $X$ is uniformly ramifying and has property (P2).
\item $X$ is almost transitive and has property (P1).
\item $X$ is almost transitive and has property (P2).
\end{enumerate}
\end{thm}
\begin{proof}
The implications (3) $\Rightarrow$ (2) and (5) $\Rightarrow$ (4) are trivial. Lemma \ref{al_trans_then_ram_unif} implies (4) $\Rightarrow$ (2) and (5) $\Rightarrow$ (3). We will now prove (1) $\Rightarrow$ (5) and (2) $\Rightarrow$ (1).\par
(1) $\Rightarrow$ (5) \quad If $X$ is quasi-isometric to $T$ then, by Lemma \ref{if_qi_then_rays...}, there cannot exist a mixed end in $\Omega X$. Thick proper ends cannot be mapped onto $\Omega T$ under $\Phi$ since thick ends cannot be described by a sequence of $n$-cuts. By Lemma \ref{if_qi_then_rays...} they also cannot be mapped onto $VT$ and therefore $\Theta_{2}X=\emptyset$.\par
Another consequence of Lemma \ref{if_qi_then_rays...} is $\Phi(\Omega_{0})\subset VT$. The only remaining ends with rays of infinite diameter are the thin proper ends. By Lemma \ref{qi_the_Phi_bij_on...} $\Phi$ is bijective on $\Phi^{-1}(\Omega T)$. Thus we have $\Phi(\Delta_{2}X)=\Omega T$. Now it is also clear that $\Phi^{-1}(VT)$ equals $\Omega_{0}X$.\par
By Lemma \ref{if_qi_then_at} $X$ is almost transitive.\par
(2) $\Rightarrow$ (1) \quad Assuming that $\Phi^{-1}(VT)=\Omega_{0}X$ we want to prove that a graph with uniform ramification is quasi-isometric under the vertex structure mapping $\phi$ to a structure tree $T=T(E)$.\par
We suppose that there exists a vertex $v$ in $VT$ having region $R(v)$ of infinite diameter. The stabilizer $\Aut_{v}^{T}(X)$ of $v$ maps the neighbours of $v$ in $T$ onto themselves. The set $L^{-1}(\Aut_{v}^{T}(X))$ of the corresponding automorphisms in $\Aut(X)$ has at most two orbits $O_{1}$ and $O_{2}$ on
$$N^{*}(v)=\{e\in E\mid o(e)=v\},$$
where $L$ is the function defined in Section \ref{act_of_Aut_on}.
For some cut $e\in O_{1}$ we now choose a finite and connected subgraph $W_{1}$ of $X$ with $VW_{1}\subset e$ so that $W_{1}$ connects all pairs of vertices in $I\theta e$ by paths of minimal length that do not leave $e$. Note that these paths are not necessarily geodesic. We define
$$\bar{W_{1}}:=\bigcup_{g\in \Aut_{e}(X)}gW_{1}.$$
$\bar{W_{1}}$ has finite diameter, because the distance of every vertex in $W_{1}$ to $\theta e$ is at most $\diam_{X}W_{1}$. For every $f\in O_{1}$ we now replace the restriction of $X$ onto $f$ by an automorphic image of $\bar{W_{1}}$. With the possibly existing orbit $O_{2}$ we proceed analogously. Thereby we obtain a connected subgraph $\bar{X}$ of $X$. We set
$$\bar{O_{i}}:=\{Vg\bar{W}_{i}\mid g\in L^{-1}(\Aut_{v}^{T}(X))\}\quad i=1,2$$
$$\mbox{and\quad }\bar{O}:=\bar{O_{1}}\cup\bar{O_{2}}.$$
Let $x$ and $y$ be two vertices in $R(v)$ and $P(x,y)$ a $\dm_{X}$-geodesic path connecting them. All parts of maximal length in $P(x,y)$ that are completely contained in one of the orbits $O_{1}$ and $O_{2}$ can be replaced by a path in $\bar{O}$ of the same length. Thus we have
$$\dm_{\bar{X}}(x,y)=\dm_{X}(x,y)$$
for all pairs of vertices $x$ and $y$ in $R(v)$.\par
The automorphism group $\Aut(\bar{X})$ generates no more than two orbits on $\bar{O}$. In order to prove the existence of a ray with infinite diameter, we now proceed analogously to the proof of Lemma \ref{al_trans_then_ram_unif}. Assuming that there is no ray with infinite diameter in $\bar{X}$, by Lemma \ref{diam_of_graphs}, there must exist some star ball $S=B(z,r)$. By Lemma \ref{ball_contains_sb_then_sb}, $S$ can be chosen so that it contains an element of both $\bar{O_{1}}$ and $\bar{O_{2}}$. There is no radius $r_{1}$ such that all sets in $\bar{O}$ which are contained in components of $\mathcal{C}_{0}(S)$ are subsets of $B(z,r_{1})$, because then, by Lemma~\ref{ball_contains_sb_then_sb}, $B(z,r_{1})$ would also be a star ball in $X$, which would be a contradiction to the uniform ramification of $X$.\par
Let $\bar{w}\in \bar{O}$ be a set of vertices which is contained in a component $C$ in $\mathcal{C}_{0}(S)$ such that its distance to $S$ is at least $2r$. Since $S$ contains an element of both orbits $\bar{O_{1}}$ and $\bar{O_{2}}$ there must exist an automorphism $g\in \Aut(\bar{X})$ such that $g (\bar{w})$ is contained in $S$. Now $g^{-1}(S)$ is completely contained in $C$. $V\bar{X}\backslash C$ is connected and therefore part of a component in $\mathcal{C}(g^{-1}(S))$. But this is impossible, because then $V\bar{X}\backslash g^{-1}(S)$ would have one component of infinite diameter, all other components would be contained in $C$ and $g^{-1}(S)$ would not be a star ball in $\bar{X}$.\par
Thus, by Lemma \ref{al_trans_then_ram_unif}, there exists a ray $L$ with infinite diameter in $\bar{X}$. Since $\dm_{\bar{X}}(x,y)=\dm_{X}(x,y)$ for all vertices $x$ and $y$ in $R(v)$, the ray $L$ also has infinite diameter in $X$. The ray $L$ has finite intersection with every cut in $N(v)^{*}$, because the intersection of these cuts with $V\bar{X}$ is finite. Hence the $\Phi$-image of the end of $L$ must be $v$. This is a contradiction to the condition $\Phi^{-1}(VT)=\Omega_{0}X$.\par
We have now proved that, assuming $\Phi^{-1}(VT)=\Omega_{0}X$, there is no vertex in $VT$ whose region has infinite diameter. So, by Theorem \ref{general_qi-criterion}, the graph $X$ is quasi-isometric to its structure tree $T$.
\end{proof}

\section{Another criterion for quasi-isometry}

The stabilizer $\Aut_{\omega}(X)$ of an $\omega\in \Omega X$ is the group of automorphisms in $\Aut(X)$ that map rays and vertices in $\omega$ onto rays and vertices in $\omega$.\par
The following theorem was proved by M\"oller for locally finite graphs with infinitely many ends where the stabilizer $\Aut_{\omega}(X)$ acts transitively. See \cite[Theorem 1]{moeller92ends2}.

\begin{thm}\label{end_stab}
Let\/ $X$ be a connected graph with a structure tree\/ $T=T(E)$. If there is an end\/ $\omega\in \Omega X$ such that the stabilizer\/ $\Aut_{\omega}(X)$ acts almost transitively on\/ $X$, then\/ $X$ is quasi-isometric to\/ $T$ by the vertex structure mapping\/ $\phi$.
\end{thm}

\begin{proof}
If $\Aut_{\omega}(X)$ acts almost transitively on $X$, then also $\Aut_{\omega}^{T}(X)$ must act almost transitively on $T$. Thus $\Phi(\omega)$ must be an end in $\Omega T$.\par
First we prove that all $\phi$-pre-images of vertices in $VT$ have finite diameter. If there is a constant $n_{0}$ such that any two vertices in $VX$ with $\dm_{X}$-distance at least $n_{0}$ have a different $\phi$-image, then $\diam_{X}\phi^{-1}(v)\le n_{0}$ for all vertices $v\in VT$. This is equivalent to the condition that for all pairs of vertices in $VX$ with $\dm_{X}$-distance at least $n_{0}$ there exists a cut in $E$ which separates them. Let $f\in ET$ be a cut containing $\omega$ and let $B(x_{0},r_{0})$ be a covering ball of $X$ with respect to $\Aut_{\omega}(X)$. We define
$$M_{f}:=\{x\in f^{*}\mid \dm_{X}(f,x)\le 4 r_{0}\}.$$
$M_{f}$ has a finite diameter. Let $y$ be a vertex in $M_{f}$ with $\dm_{X}(f,y)=2r_{0}$. The ball $B(y,r_{0})$ contains a vertex $y_{0}$ of the $x_{0}$-orbit with respect to $\Aut_{\omega}(X)$. As $B(y_{0},r_{0})$ is also a subset of $M_{f}$, we have
$$B(y_{0},r_{0})\subset B(y,2 r_{0})\subset M_{f}.$$
We define $n_{0}:=2 \diam_{X} M_{f}$ and choose two arbitrary vertices $x_{1}$ and $x_{2}$ with distance larger then $n_{0}$. Since $\bigcup \Aut_{\omega}(X) M_{f}$ is the whole set of vertices $VX$, there is a cut $e_{1}$ containing $\omega$ for which $x_{1}\in M_{e_{1}}$. If $x_{2}$ is an element of $e_{1}$ there is nothing more to prove because then $e_{1}$ is the desired cut which separates $x_{1}$ and $x_{2}$. So we suppose that $x_{2}$ is an element of $e_{1}^{*}$. Let $e_{2}$ be a cut that contains $\omega$ and for which $x_{2}\in M_{e_{2}}$. The edge-boundaries $\delta e_{1}$ and $\delta e_{2}$ are disjoint and $M_{e_{2}}$ is a subset of $e_{1}^{*}$. Thus $e_{1}\cup M_{e_{1}}$ must be a subset the component of $e_{2}$ which contains $\omega$. Thus $x_{1}\in e_{2}$ and $x_{2}\in e_{2}^{*}$.\par
To prove the theorem we have to show that the region of any vertex $v\in VT$ has finite diameter. By Lemma \ref{diam_phi^-1_and_B} we just have to deal with the case $\phi^{-1}(v)=\emptyset$. Let $L:=\{w_{0}=v,w_{1},w_{2},\dots\}$ be the ray which starts at $v$ and lies in $\Phi(\omega)$. We furthermore define
$$U(v):=\{u\in VT\mid u\sim v \mbox{\quad and\quad } u\ne w_{1}\}.$$
The $\phi$-pre-image of any vertex $v_{0}$ in $U(v)$ is non-empty and has a finite diameter. Thus it is contained in some covering ball $B$ with respect to $\Aut_{\omega}(X)$. Let $M_{0}$ be the set of all $\phi$-pre-images of a vertices in $VT$ which have a non-empty intersection with $B$. The diameter of $M_{0}$ is finite. By Lemma \ref{general_Q1} this also holds for $\phi(M_{0})$. Let $n$ be the smallest index such that no vertex in $\{w_{n},w_{n+1},w_{n+2},\dots\}$ lies in $\phi(M_{0})$. We define
$$d:=\diam_{T}\phi(M_{0})\mbox{\quad and}$$
$$A:=\{w_{n},w_{n+1},\dots,w_{n+2d}\}.$$
For every $u\in U(v)$ there exists an automorphism $g_{u}\in \Aut^{T}_{\omega}(X)$ with $g_{u}(u)\in \phi(M_{0})$. Since such an automorphism $g_{u}$ must fix $\Phi(\omega)$ and  $u$ to a vertex in $VT$ which has a distance to $u$ that is at most $d$, it causes a translation on $L$ of maximal length $d$ and therefore it must map $w_{n+d}$ onto a vertex in $A$. Hence
$$\max \{\dm_{X}(\phi^{-1}(u),\phi^{-1}(w_{n+d}))\mid u\in U(v)\}\le
\max \{\dm_{X}(\phi^{-1}(v_{0}),\phi^{-1}(w))\mid w\in A\}<\infty.$$
This implies
$$\diam_{X}\phi^{-1}(w_{1})\cup \bigcup_{u\in U(v)} \phi^{-1}(u))<\infty.$$
Since
$$R(v)\subset \phi^{-1}(w_{1})\cup \bigcup_{u\in U(v)} \phi^{-1}(u)=\phi^{-1}(B(v,1))$$
we finally have
$$\diam_{X}R(v)<\infty.$$
\end{proof}

\begin{exmp}\par
\mbox{}
\begin{enumerate}
\item For a semi-regular tree $T$ there is, up to isomorphism, a unique pair $e$ and $e^{*}$ of structure cuts. The corresponding structure tree is isomorphic to $T$.
\item Let $T$ be the tree in 1. By adding a graph of finite diameter to all the vertices in one of the bipartite blocks of $T$ we obtain a graph $X$ which is quasi-isometric to $T$.
\item Let $X$ be the Cayley graph of the free product $\mathbb{Z} * \mathbb{Z}^{2}=\left<a,b,c\mid bc=cb\right>$ with generating system $\{a^{\pm 1}, b^{\pm 1}, c^{\pm 1}\}$. This graph is 6-regular and transitive. By removing edges that correspond to the generating elements $a^{\pm 1}$ we obtain pairs of structure cuts. The structure tree $T$ is regular of countably infinite degree. $\Phi$ maps thick ends onto vertices and thin ends onto thin ends. $T$ and $X$ are not quasi-isometric.\par
Taking $\{a^{\pm 1}\}\cup \mathbb{Z}^{2}$ as a generating system we obtain the Cayley graph $\bar{X}$. The ends in the copies of $\mathbb{Z}^{2}$ now are thick point ends. Again removing the edges that correspond to the generating elements $a^{\pm 1}$ we obtain a structure tree which is isomorphic to $T$. $\bar{X}$ and $T$ are quasi-isometric, because $\Phi$ maps only point ends onto vertices in $VT$.
\end{enumerate}
\end{exmp}

In \cite{trofimov85groups} Trofimov gave an example of a graph whose automorphism group fixes an end and acts transitively on the set of vertices.\newline

\begin{appendix}
This article is based on Chapter 3 of the author's masters thesis \cite{kroen98topologische} at the University of Salzburg under supervision of Prof.\ W.\ Woess, and the author wants to thank him for many useful suggestions. The main part of this thesis was written during a stay at Milan supported by the Italian Ministry of Foreign Affairs.
\end{appendix}
\newline

\end{document}